\documentclass[12pt,twoside]{article}
\newcommand{\byd}{\stackrel{\rm d}{\rightarrow}}
\newcommand{\byp}{\stackrel{\rm P}{\longrightarrow}}
\newcommand{\lip}{\mathop{\rm l.i.p.}\limits_{n\rightarrow\infty}}
\newcommand{\h}{{\bf H}}
\newcommand{\li}{{\rm Lip}}

\newcommand{\R}{{\rm R}}
\newcommand{\pr}{{\rm P}}

\begin{document}
\thispagestyle{plain}
\newcounter{N}

\noindent{\small Theory of Stochastic Processes \\ Vol.12 (28),
no.3-4, 2006, pp.*-*} \vspace{2cm}\addtolength{\topmargin}{-0.7cm}

\markboth{\hfil {\small\rm DMYTRO IVANENKO
\hbox{\vrule\hspace{2.5cm}}
}\hfil}{\hfil{\small\rm\hbox{\vrule\hspace{2cm}} ASYMPTOTICALLY
OPTIMAL ESTIMATOR}\hfil}

\begin{center}
DMYTRO IVANENKO
\end{center}

\begin{center}
{\large \bf ASYMPTOTICALLY OPTIMAL ESTIMATOR OF THE PARAMETER OF
SEMI-LINEAR AUTOREGRESSION}
\footnote{%Invited lecture.

}
\end{center}

\begin{quote}
{\small  The difference equations
$\xi_{k}=af(\xi_{k-1})+\varepsilon_{k}$, where $(\varepsilon_k)$
is a square integrable difference martingale, and the differential
equation ${\rm d}\xi=-af(\xi){\rm d}t+{\rm d}\eta$, where $\eta$
is a square integrable martingale, are considered. A family of
estimators depending, besides the sample size $n$ (or the
observation period, if time is continuous) on some random
Lipschitz functions is constructed. Asymptotic optimality of this
estimators is investigated.}
\end{quote}

\begin{center}
{\sc 1. Introduction}
\end{center}

\begin{center}
{\it Discrete time}
\end{center}

We consider the difference equation
\begin{equation}\label{def1}
 \xi_{k}=af(\xi_{k-1})+\epsilon_{k} , \quad
k\in\mathrm{N},
\end{equation}
where $\xi_0$ is a prescribed random variable, $f$ is a prescribed
nonrandom function,  $a$ is an unknown scalar parameter and
$(\epsilon_k)$ is a square integrable difference martingale with
respect to some flow $\left( {\rm F}_{k}, k\in\mathrm{Z}_+\right)$
of $\sigma$-algebras such that the random variable $\xi_0$ is
${\rm F}_0$-measurable. In the detailed form, the assumption about
$(\epsilon_k)$ means that for any $k$ $\epsilon_k$ is ${\rm
F}_k$-measurable,
\begin{equation}\label{cond1}
 {\rm E}\epsilon_k^2<\infty
\end{equation} and
\begin{equation}\label{cond2}
 {\rm E}(\epsilon_k|{\rm F}_{k-1})=0.
\end{equation}

The word "semi-linear"~in the title means that the right-hand
sides of (\ref{def1}) depend linearly on $a$ but not on $\xi$.

We use the notation: ${\rm l.i.p.}$ -- limit in probability;
$\byd$ -- the weak convergence of finite-dimensional distributions
of random functions, in particular convergence in distribution of
random variables.

Let for each $k\in{\rm Z}_+$ $h_k=h_k(\omega, x)$ be an ${\rm
F}_{k-1}\otimes{\rm B}$-measurable function (that is the sequence
$(h_k)$ be predictable) such that
$$
 {\rm E}\left[\left(|\xi_{k+1}|+|af(\xi_{k+1})|\right)|h_k(\xi_{k})|\right]+{\rm E}|h_k(\xi_{k})|<\infty.
$$
Then from (\ref{def1}) -- (\ref{cond2}) we have ${\rm
E}\left(\xi_{k+1}-af(\xi_{k})\right) h_k(\xi_{k})=0$, whence
$$
 a=\left({\rm E}\xi_{k+1} h_k(\xi_{k})\right)\left({\rm E} f(\xi_{k})
 h_k(\xi_{k})\right)^{-1}
$$
provided $\left({\rm E} f(\xi_{k})
 h_k(\xi_{k})\right)\neq0$. This prompts the estimator
\begin{equation}\label{estim}
  \check{a}_n=\left(\sum_{k=0}^{n-1}\xi_{k+1}
  h_{k}(\xi_{k})\right)\left(\sum_{k=0}^{n-1}
  f(\xi_{k})
 h_k(\xi_{k})\right)^{-1},
\end{equation}
coinciding with the LSE if $h_k(x)=f(x)$ for all $k$.

\begin{center}
{\it Continuous time}
\end{center}

We consider the differential equation
\begin{equation}\label{def2}
 {\rm d}\xi(t)=-af(\xi{(t)}){\rm d}t+{\rm d}\eta(t), \quad
 t\in{\mathrm R},
\end{equation}
where $\eta(t)$  is a local square integrable martingale w.r.t. a
flow ${(\rm F}(t))$ such that the random variable $\xi(0)$ is
${\rm F}(0)$-measurable.

Let $h(t,x)$ be a predictable random function such that for all
$t\in{\mathrm R}_+$
$$
 {\rm E}\left[\left(|\xi{(t)}|+|af(\xi{(t)})|\right)|h(t,\xi{(t)})|\right]+{\rm E}|h(t,\xi{(t)})|<\infty.
$$
Let us multiply (\ref{def2}) on $h(t, \xi(t))$ and integrate from
$0$ to $T$. The same rationale as in the discrete case yields the
estimator
\begin{equation}\label{estc}
 \check{a}_T=-\left(\int_0^T h(t,\xi(t)){\rm d}\xi\right)\left(\int_0^T f(\xi(t))h(t,\xi(t)){\rm
 d}t\right)^{-1},
\end{equation}
coinciding with the LSE if $h(t,x)=f(x)$.
\smallskip

Asymptotic normality of $\sqrt{n}\left(\check{A}_{n}-A\right)$,
where $\check{A}_{n}$ is the LSE of a matrix parameter $A$, was
proved in [1] under the assumptions of ergodicity and stationarity
of $(\xi_n)$. Convergence in distribution of this normalized
deviation was proved in [2] with the use of stochastic calculus.
Ergodicity and even stationarity of $(\epsilon_k)$ was not assumed
in [2], so the limiting distribution could be other than normal.

The goal of the article is to match a sequence $(h_k)$ (if time is
discrete) or a function $h(t,\cdot)$ (if time is continuous) so
that to minimize the value of some random functional $V_n$ which,
as we shall see in Section 3, is asymptotical close in
distribution to some numeral characteristic of the estimator (in
case the latter is asymptotically normal this characteristic
coincides with the variance).

\begin{center}
{\sc 2. The main results}
\end{center}

\begin{center}
{\it Discrete time}
\end{center}

Denote $\sigma_k^2=\rm{E[\epsilon_k^2|F_{k-1}]}$,
$\mu_k=h_k(\xi_k)$. Let $\li(C)$ denote the class of functions
satisfying the Lipschitz condition with some constant $C$ and
equal to zero at the origin, $\li={\bigcup}_{C>0}\li(C)$, and let
$\h(C)$ denote the class of all predictable random functions on
${\rm Z_+\times R}$ (discrete time) or ${\rm R_+\times R}$
(continuous time) whose realizations $h_k(\cdot)$ (respectively
$h(t, \cdot)$) belong, as functions of $x$, to $\li(C)$,
$\h={\bigcup}_{C>0}\h(C)$. Predictability means ${\rm P\otimes
B}$-measurability in $(\omega, t, x)$ (the $\sigma$-algebra $\pr$
is defined in [4, p. 28], [6, p.~13]).

We are seeking for $(\widetilde{h}_k)\in\h$ minimizing the
functional
\begin{equation}\label{disp}
    V_n(h_0,\ldots, h_{n-1})=\frac{\frac{1}{n}\sum_{k=0}^{n-1}\sigma_{k+1}^2\mu_k^2}
{\left(\frac{1}{n}\sum_{k=0}^{n-1}f(\xi_{k})\mu_k\right)^2}.
\end{equation}
\smallskip

\noindent{\bf Theorem 1.} {\it Let
\begin{equation}\label{prob}
V_n(\widetilde{h}_0,\ldots,
\widetilde{h}_{n-1})=\min_{h_0,\ldots,h_{n-1}\in
{\h}}V_n(h_0,\ldots, h_{n-1}).
\end{equation}
Then
\begin{equation}\label{equation}
\sigma_{k+1}^2\widetilde{\mu}_k\sum_{i=0}^{n-1}f(\xi_{i})\widetilde{\mu}_i
=f(\xi_k)\sum_{i=0}^{n-1}\sigma_{i+1}^2\widetilde{\mu}_i^2, \qquad
k=\overline{0,n-1}.
\end{equation}}

\noindent{\it Proof.} To obtain the necessary conditions for
extremum of the functional $V_n$ (\ref{equation}) we will vary [3]
just one of functions $h_k,~k=\overline{0,n-1}$, leaving the other
functions without changes. Thus regarding $V_n(h_0,\ldots,
h_{n-1})$ as a functional depending on only one function
$V_n(h_0,\ldots, h_{n-1})=\widetilde{V}_n(h_k)$.

Let's choose some scalar function $g\in \h$ and denote
$g_\lambda(x)=\widetilde{h}_k(x)+\lambda(g(x)-\widetilde{h}_k(x))$,
$v(\lambda)=\widetilde{V}_n(g_\lambda)$.

Obviously, $g_\lambda\in\h$ so the minimum of $v(\lambda)$ is
attained at zero and therefore
\begin{equation}\label{pro1}
    v'(0)=0.
\end{equation}
The expression for the left-hand side is
$$
v'(0)=
\frac{2n(g(\xi_k)-\mu_k)\left(\sigma_{k+1}^2\widetilde{\mu}_k(\sum_{i=0}^{n-1}f(\xi_{i})\widetilde{\mu}_i-
 f(\xi_{k})\sum_{i=0}^{n-1}\sigma_{i+1}^2\widetilde{\mu}_i^2\right)}
{\left(\sum_{i=0}^{n-1}f(\xi_{i})\widetilde{\mu}_i \right)^3}.
$$
Hence in view of (\ref{pro1}) we obtain the $i$ th equation of
system (\ref{equation}).

It remains to apply this argument to each function $h_k,~
k=\overline{0,n-1}$.
\smallskip

\noindent{\bf Remark.} The Lipschitz condition was not used in the
proof. It will be required in Section 3.

\noindent{\bf Corollary 1.} Let $f\in {\li(C)}$ and there exist a
constant $q> 0$ such that $\sigma_k^2\geq q$ for all $k$. Then
${h}_i(x)=f(x)/\sigma_{i+1}^2$, $i=\overline{0,n-1}$, is a
solution to the problem (\ref{prob}).

\begin{center}
{\it Continuous time}
\end{center}

Let $m$ denote the quadratic characteristic of $\eta$.

We shall match $\widetilde{h}=\widetilde{h}(\omega, t, x)$ from
${\h(C)}$ ($C$ is independent of $t$) so that to minimize the
value of the functional
\begin{equation}\label{disp1}
    V_T(h)=\frac{\frac{1}{T}\int_{0}^{T}h(t,\xi{(t)})^2{\rm d}m(t)}
{\left(\frac{1}{T}\int_{0}^{T}f(\xi{(t)})h(t,\xi{(t)}){\rm
d}t\right)^2}.
\end{equation}

\noindent{\bf Theorem 2.} {\it Let %$\widetilde{h}$ be a solution
%to the problem
\begin{equation}\label{prob1}
V_T(\widetilde{h})=\min_{h\in \h}V_T(h).
\end{equation}
Then for all $g\in {\h}$
\begin{equation}\label{equation1}
\begin{array}{l} \int_{0}^T\widetilde{h}(t,\xi(t))g(t,\xi(t)){\rm
d}m(t)\int_{0}^{T}f(\xi{(t)})\widetilde{h}(t,\xi{(t)}){\rm d}t =
\\ \qquad\qquad\qquad\qquad\qquad\qquad\int_0^Tf(\xi(t))g(t,\xi(t)){\rm
d}t\int_{0}^{T}\widetilde{h}(t,\xi{(t)})^2{\rm d}m(t).\end{array}
\end{equation}}

\noindent{\it Proof.} Let's choose some scalar function $g\in \h$
and denote $g_\lambda(t,x)=\widetilde{h}(t,x)+\lambda g(t,x)$,
$v(\lambda)={V}_T(g_\lambda)$.

Obviously $g_\lambda(t, \cdot)\in\h$ so the minimum of
$v(\lambda)$ is attained in zero and therefore
\begin{equation}\label{pr2}
    v'(0)=0.
\end{equation}
The expression for the left-hand side is
\newline $v'(0)=2T\left(\int_{0}^{T}f(\xi{(t)})\widetilde{h}(t,\xi{(t)}){\rm
d}t\right)^{-3}\times$
$$\left(\int_{0}^{T}f(\xi(t))\widetilde{h}(t,\xi{(t)}){\rm d}t\int_{0}^T\widetilde{h}(t,\xi(t))g(t,\xi(t)){\rm
d}m(t)-\right.$$
$$\quad\qquad\qquad\qquad\qquad\qquad\qquad\left.\int_0^Tf(\xi(t))g(t,\xi(t)){\rm
d}t\int_{0}^{T}\widetilde{h}(t,\xi{(t)})^2{\rm d}m(t)\right).
$$
Hence in view of (\ref{pr2}) we come to (\ref{equation1}).
\smallskip

\noindent{\bf Corollary 2.} Let $f\in {\li(C)}$, $m$ be absolutely
continuous w.r.t. the Lebesgue measure and there exist a constant
$q> 0$ such that for all $t$ $\dot{m}\geq q$. Then
${h}(t,x)=f(x)/\dot{m}$ is a solution to the problem
(\ref{prob1}).

\begin{center}
{\sc 3. An illustration }
\end{center}

Denote ${\rm E}^0={\rm E}(\cdots|{\rm F}_0)$,
 $Q_n=\frac{1}{n}\sum_{k=0}^{n-1}f(\xi_k)\mu_k$,
 $G_n=\frac{1}{n}\sum_{k=1}^{n}\sigma_k^2\mu_{k-1}^{2}$.

We denote ${\rm E}^0={\rm E}(\cdots|{\rm F}_0)$ and introduce the
conditions
\newline
{\bf CP1.} For any $r\in {\rm N}$ and any uniformly bounded
sequence $(\alpha_k)$ of ${\mathrm R}$-valued Borel functions on
${\mathrm R}^{r}$
$$
 \frac{1}{n}\sum_{k=r}^{n-1}\left(\alpha_k(\epsilon_{k-r+1},\ldots ,\epsilon_k)-
 {\rm E}^0\alpha_k(\epsilon_{k-r+1},\ldots
 ,\epsilon_k)\right)\byp0,
$$
$$
 \frac{1}{n}\sum_{k=r}^{n-1}\left(\sigma_k^2\alpha_k(\epsilon_{k-r+1},\ldots ,\epsilon_k)-
 {\rm E}^0\sigma_k^2\alpha_k(\epsilon_{k-r+1},\ldots
 ,\epsilon_k)\right)\byp0.
$$
{\bf CP2.} For such $r$ and $(\alpha_k)$ the sequences
$$
 \left(\frac{1}{n}\sum_{k=r}^{n-1}{\rm E}^0\alpha_k(\epsilon_{k-r+1},\ldots ,\epsilon_{k}),
 \quad n=r+1,\ldots\right),
$$
$$
 \left(\frac{1}{n}\sum_{k=r}^{n-1}{\rm E}^0\sigma_k^2\alpha_k(\epsilon_{k-r+1},\ldots ,\epsilon_{k}),
 \quad n=r+1,\ldots\right)
$$
converge in probability.

Denote $f_0(x)=x$ and, for $r\geq1$,
$$
  f_r(x_0,\ldots ,x_r)=af(f_{r-1}(x_0,\ldots ,x_{r-1}))+x_r.
$$
Then
$$
  \xi_k=f_r(\xi_{k-r},\epsilon_{k-r+1},\ldots
  ,\epsilon_k),\quad r<k.
$$
\smallskip

\noindent{\bf Lemma 1.} {\it  Let conditions (\ref{cond1}),
(\ref{cond2}), {\bf CP1} and {\bf CP2} be fulfilled. Suppose also
that
\begin{equation}\label{condeps}
  \lim_{N\rightarrow\infty}\overline{\lim_{n\rightarrow\infty}}
  \frac{1}{n}\sum_{k=1}^{n}{\rm E}\epsilon_k^2 I\{|\epsilon_k|>N\}=0
\end{equation}
and there exist an ${\rm F}_0$-measurable random variable
$\upsilon$ such that for all $k$
\begin{equation}\label{condeps1}
 \sigma_k^2\leq \upsilon.
\end{equation}
and positive numbers $C,~C_1$ such that
\begin{equation}\label{conda}
    |a|C< 1,
\end{equation}
$f\in {\li(C)}$, $(h_k)\in\h(C_1)$. Then
\begin{equation}\label{pr1}
 (G_n, Q_n)\byd(G, Q).
\end{equation}}
\smallskip

\noindent{\it Proof.} Denote
$\xi_k^r=f_r(0,\epsilon_{k-r+1},\ldots ,\epsilon_k)$,
$\mu_k^r=h_k(\xi_k^r)$,
$Q_n^r=\frac{1}{n}\sum_{k=r}^{n-1}f(\xi_k^r)\mu_k^r$,\newline
$G_n^r=\frac{1}{n}\sum_{k=r}^{n}\sigma_k^2(\mu_{k-1}^r)^{2}$. We
claim that conditions (\ref{cond1}), (\ref{cond2}),
(\ref{condeps}), (\ref{condeps1}), (\ref{conda}) and the relation
\begin{equation}\label{pro3}
 (Q_n^r, G_n^r)\byd (Q^r, G^r)\quad{\rm as}\quad
 n\rightarrow\infty
\end{equation}
imply (\ref{pr1}).

Let $X_r$ denote $(x_1,\ldots, x_r)\in{\rm R}^r$. Then under the
assumptions on $f$ and $h_k$ for any $N>0$
$$
 \lim_{r\rightarrow\infty}\sup_{|x|\leq N, X_r\in\R^{r}}|f_{r}(x, X_{r})-f_{r}(0,
 X_{r})|=0,
$$
whence with probability 1 for any $k$
\begin{equation}\label{pr5}
 \lim_{r\rightarrow\infty}\sup_{|x|\leq N, X_r\in\R^{r}}|f(f_{r}(x, X_{r})) h_k(f_{r}(x,
 X_{r}))-f(f_{r}(0, X_{r}))h_k(f_{r}(0,
 X_{r}))|=0,
\end{equation}
$$
 \lim_{r\rightarrow\infty}\sup_{|x|\leq N, X_r\in\R^{r}}|h_k(f_{r}(x, X_{r}))^{2}-h_k(f_{r}(0,
 X_{r}))^{2}|=0.
$$
These relations was proved in [5].

Let us prove that from conditions (\ref{cond1}), (\ref{cond2}),
(\ref{condeps}), (\ref{condeps1}) and (\ref{conda}) it follows
that almost surely
\begin{equation}\label{pr4}
  \lim_{r\rightarrow\infty}\overline{
  \lim_{n\rightarrow\infty}}{\rm E}^0|Q_n-Q_n^r|=0,\qquad
  \lim_{r\rightarrow\infty}\overline{
  \lim_{n\rightarrow\infty}}{\rm E}^0|G_n-G_n^r|=0.
\end{equation}

By (\ref{pr5}) for any $N>0$
\begin{equation}\label{pr6}
  \lim_{r\rightarrow\infty}\overline{
  \lim_{n\rightarrow\infty}}\frac{1}{n}\sum_{k=r}^{n-1}{\rm E}|f(\xi_k)
  \otimes\mu_k-f(\xi_k^r)\mu_k^r|I\{|\xi_k|\leq
  N\}=0.
\end{equation}
Denote $ \chi_{k}^{N}=I\{|\xi_{k}|>N\},$
  $I_{k}^{N}=I\{|\epsilon_{k}|>(1-C)N\},$
  $ b_{k}^{N}={\rm E}^{0}|\xi_{k}|^{2}\chi_{k}^{N}.$
Due to (\ref{conda}) and because of $(h_k)\in\h(C_1)$
$$
{\rm E}^0|f(\xi_k)\mu_k|\chi_k^N\leq CC_1b_k^N,
$$
Hence and from (\ref{cond1}), (\ref{cond2}),
(\ref{condeps})--(\ref{conda}) we get by Corollary 1 [5]
\begin{equation}\label{pr7}
 \lim_{N\rightarrow\infty}\overline{
  \lim_{n\rightarrow\infty}}\frac{1}{n}\sum_{k=0}^{n-1}{\rm E}^0
 |f(\xi_k)\mu_k|\chi_k^N=0.
\end{equation}
Further, for $k\geq r$,
$$
 {\rm E}^0|f(\xi_k^r)\mu_k^r|={\rm E}^0|f(f_r(0,\epsilon_{k-r+1},\ldots ,\epsilon_k))|
 |h_k(f_r(0,\epsilon_{k-r+1},\ldots ,\epsilon_k))|,
$$
whence
\begin{equation}\label{pr8}
 {\rm E}|f(\xi_k^r)\mu_k^r|\chi_k^N\leq
 CC_1{\rm E}\left(\sum_{i=0}^{r-1}C^i|\epsilon_{k-i}|\right)^2\chi_k^N.
\end{equation}
Writing the Cauchy -- Bunyakovsky inequality
$$
 \left(\sum_{i=0}^{r-1}C^i|\epsilon_{k-i}|\right)^2\leq \sum_{j=0}^{r-1}C^j
 \sum_{i=0}^{r-1}C^i|\epsilon_{k-i}|^2,
$$
we get for an arbitrary $L>0$ \newline${\rm
E}\left(\sum_{i=0}^{r-1}C^i|\epsilon_{k-i}|\right)^2\chi_k^N\leq$
 \begin{equation}\label{pr9} (1-C)^{-1}
  \left({\rm E}\sum_{i=0}^{r-1}C^i\epsilon_{k-i}^2I\{|\epsilon_{k-i}|>L\}
  +L^2\pr\{|\xi_k|>N\}\sum_{i=0}^{r-1}C^i\right).
\end{equation}
In view of (\ref{cond1}), (\ref{cond2}) Lemma 1 [5] together with
(\ref{conda}) and (\ref{condeps}) implies that
\begin{equation}\label{pr10}
 \lim_{N\rightarrow\infty}\overline{
  \lim_{n\rightarrow\infty}}\frac{1}{n}\sum_{k=0}^{n}\pr\{|\xi_k|>N\}=0.
\end{equation}
Obviously, for arbitrary nonnegative numbers
$u_0,\ldots ,u_{r-1},v_1,\ldots , v_{n-1}$
$$
 \sum_{k=r}^{n-1}\sum_{i=0}^{r-1}u_iv_{k-i}\leq\sum_{i=0}^{r-1}u_i\sum_{j=1}^{n-1}v_j,
$$
so conditions (\ref{conda}) and (\ref{condeps}) imply that
$$
 \lim_{L\rightarrow\infty}\sup_r\overline{
  \lim_{n\rightarrow\infty}}\frac{1}{n}\sum_{k=r}^{n-1}
  {\rm E}\sum_{i=0}^{r-1}C^i\epsilon_{k-i}^2I\{|\epsilon_{k-i}|>L\}=0,
$$
whence in view of (\ref{pr8}) -- (\ref{pr10})
$$
 \lim_{N\rightarrow\infty}\sup_r\overline{
  \lim_{n\rightarrow\infty}}\frac{1}{n}\sum_{k=r}^{n-1}{\rm E}|f(\xi_k^r)\mu_k^r|\chi_k^N=0.
$$
Combining this with (\ref{pr6}) and (\ref{pr7}), we arrive at the
first relation of (\ref{pr4}).

The proof of the second relation of (\ref{pr4}) is similar.

The details can be found in [5].

From (\ref{pro3}), and (\ref{pr4}) we obtain that the sequence
$((Q^r, G^r), r\in{\rm N})$ converges in distribution to some
limit $(Q, G)$ and relation (\ref{pr1}) holds.

Let us check (\ref{pro3}). Condition {\bf CP1} implies that
$$
  \lim_{r\rightarrow\infty}\overline{
  \lim_{n\rightarrow\infty}}{\rm E}^0|Q_n^r-{\rm E}^0Q_n^r|=0,
\qquad
  \lim_{r\rightarrow\infty}\overline{
  \lim_{n\rightarrow\infty}}{\rm E}^0|G_n^r-{\rm E}^0G_n^r|=0.
$$
It remains to note that under condition {\bf CP2} for any $r\in
{\rm N}$ the sequences $({\rm E}^0G_n^r)$ and $({\rm E}^0Q_n^r)$
converge in probability.
\smallskip

By construction $V_n(h_0,\ldots, h_{n-1})=G_nQ_n^{-2}$. The value
$Q_n=0$ is excluded by the choice of the tuple $(h_0,\ldots ,
h_{n-1})$ minimizing $V_n$.
\smallskip

\noindent{\bf Corollary 3.} Let the conditions of Lemma 1 be
fulfilled and $Q\neq0$ a.s. Then $V_n\byd V$, where $V=GQ^{-2}$.
\smallskip

Having in mind the use of stochastic analysis, we introduce the
processes $\check{a}_n(t)=\check{a}_{[nt]}$ and the flows ${\rm
F}_{n}(t)={\rm F}_{[nt]}$ with continuous time.
\smallskip

\noindent{\bf Theorem 3.} {\it  Let conditions of Lemma 1 be
fulfilled. Then
$\sqrt{n}\left(\check{a}_{n}(\cdot)-a\right)\byd\beta(\cdot)$,
where $\beta$ is a continuous local martingale with quadratic
characteristic
\begin{equation}\label{var}
 \langle\beta\rangle(t)=tV,
\end{equation}}
and initial value 0.
\smallskip

\noindent{\it Proof.} Denote $
 Y_n(t)=\frac{1}{\sqrt{n}}\sum_{k=1}^{[nt]}\epsilon_{k}\mu_{k-1}$.
Then because of (\ref{estim})
\begin{equation}\label{stat}
 \sqrt{n}(\check{a}_n(t)-a)=Y_n(t)Q_n^{-1}.
\end{equation}

By construction and conditions (\ref{cond1}), (\ref{cond2}),
(\ref{conda}) $Y_n$ is a locally square integrable martingale with
quadratic characteristic
$$
 \langle Y_n\rangle(t)=n^{-1}[n
t]G_{[nt]}.
$$

It was proved in [5]  that under conditions (\ref{cond1}),
(\ref{cond2}), (\ref{condeps}), (\ref{condeps1}), (\ref{conda})
and (\ref{pr1}) $\sqrt{n}(\check{a}_n(\cdot)-a)\byd
Y(\cdot)Q^{-1}$, where $Y$ is a continuous local martingale w.r.t.
some flow $({\rm F}(t), t\in{\mathrm R}_+)$ such that $\langle
Y\rangle(t)=t G$ and the random variable $Q$ is ${\rm
F}(0)$-measurable (and so does $G$, which can be seen from the
expression for $\langle Y\rangle$). In view of Lemma 1 it remains
to note that $V_n=\langle Y_n \rangle(1)Q_n^{-2}$ and $V=\langle Y
\rangle(1)Q^{-2}$.
\smallskip

\noindent{\bf Remark.} This theorem explains the form of
functional (\ref{disp}). In the most general case (without
conditions {\bf CP1} and {\bf CP2}) the denominator (\ref{stat})
in limit is an ${\rm F}(0)$-measurable random variable, and the
numerator tends to quadratic characteristic at the point $t=1$ of
the continuous local martingale $Y$. Thus, the numerator
(\ref{disp}) is the quadratic characteristic at $t=1$ of the
pre-limit martingale $Y_n$, and the denominator satisfies the law
of large numbers. Minimizing pre-limit variance in
$(h_k)\in\h(C_1)$, we lessen the value of limited variance of the
normalized deviation of estimator (\ref{estim}).
\smallskip

Let further ${h}_k(x)=f(x)/\sigma_{k+1}^2$. Recall that
$(h_k,~k=\overline{0,n-1})$ is a solution to the problem
(\ref{prob}). For such $h_k$ we have

\noindent{\bf Corollary 4.} Let the conditions of Corollary 1 and
Theorem 3 be fulfilled. Then

$$
 V=\left(\lim_{r\rightarrow\infty}\lip\frac{1}{n}\sum_{k=r}^{n-1}{\rm
E}^0\frac{f(\xi_{k}^r)^2}{\sigma_{k+1}^2}\right)^{-1}.
$$

\noindent{\it Proof.} Obviously $V_n=Q_n^{-1}$. By Lemma 1
$Q_n\byd Q$, where $Q=\lim\limits_{r\rightarrow\infty}\lip{\rm
E}^0 Q_n^r$. To complete the proof it remains to note that
$Q_n^r=\frac{1}{n}\sum_{k=r}^{n-1}{\rm
E}^0\frac{f(\xi_{k}^r)^2}{\sigma_{k+1}^2}$.

\begin{center}
{\sc 4. An example}
\end{center}
\smallskip

Suppose that $f\in {\li(C)}$, $h_k\in\h(C_1)$ condition
(\ref{conda}) be fulfilled. Let also
$\epsilon_n=\gamma_nb_n(\xi_{n-1})$, where $(\gamma_n)$ be a
sequence of independent random variables with zero mean and
variances $\varsigma^2_n$, $|\gamma_k|\leq C_2$, $b_n\in\h(C_3)$
and $C+C_2C_3<1$. Let also ${\rm E}\xi_0^2<\infty$

For ${\rm F}_k$ we take the $\sigma$-algebra generated by
$\xi_0;\gamma_1,\ldots,\gamma_{k}$.\newline Then
$\sigma^2_{k}=\varsigma_k^2b_k(\xi_{k-1})^2$ and $(\epsilon_n)$
satisfies (\ref{cond1}), (\ref{cond2}).

Denote further
$$ \widehat{f}_r(x_0,\ldots
,x_r)=af(\widehat{f}_{r-1}(x_0,\ldots ,x_{r-1}))+x_r
 b_r(\widehat{f}_{r-1}(x_0,\ldots ,x_{r-1})),
$$
$$\widehat{\xi}_k^r=\widehat{f}_r(0,\gamma_{k-r+1},\ldots
,\gamma_k),\quad \widehat{\mu}_k^r=h_k(\widehat{\xi}_k^r),\quad
\widehat{Q}_n^r=\frac{1}{n}\sum_{k=r}^{n-1}f(\widehat{\xi}_k^r)\widehat{\mu}_k^r,$$
$\widehat{G}_n^r=\frac{1}{n}\sum_{k=r}^{n-1}\varsigma_{k+1}^2b_{k+1}(\widehat{\xi}_k^r)^2(\widehat{\mu}_{k}^r)^{2}$.
Similarly to the proof of Lemma 1 we obtain
$$
  \lim_{r\rightarrow\infty}\overline{
  \lim_{n\rightarrow\infty}}{\rm
  E}^0|G_n-\widehat{G}_n^r|=0,\qquad\lim_{r\rightarrow\infty}\overline{
  \lim_{n\rightarrow\infty}}{\rm E}^0|Q_n-\widehat{Q}_n^r|=0.
$$

Items in $\widehat{G}_n^r$ and $\widehat{Q}_n^r$ depends on
$\gamma_{k-r+1},\ldots ,\gamma_k$ then they satisfy the law of
large numbers in Bernstein's form.

If besides $\epsilon_n$ satisfies {\bf CP2} and $Q\neq0$ then
Theorem 3 asserts (\ref{var}). If herein
$\frac{f(x)}{\varsigma_k^2b_k(x)^2}\in{\li}$ then
$\widetilde{h}_k(x)=\frac{f(x)}{\varsigma_k^2b_k(x)^2}$ is a
solution to the problem (\ref{prob}) and
$$
 V=\lim_{r\rightarrow\infty}\lip\left(\frac{1}{n}\sum_{k=r}^{n-1}{\rm
 E}^0\frac{f(\widetilde{\xi}_k^r)^2}{\varsigma_{k+1}^2b_{k+1}(\widetilde{\xi}_k^r)^2}\right)^{-1}.
$$

\smallskip

\noindent{\bf Example.} Let $b_n=b$, $h_n=h$ and $\gamma_n$ be
i.i.d. random variables. In view of expressions for
$\widehat{Q}_n^r$ and $\widehat{G}_n^r$ we may confine ourselves
with the case $\alpha_k=\alpha$.

By the Stone -- Weierstrass theorem for $\sigma$-compact spaces
[7, p. 317] $\alpha$ can be uniformly on compacta approximated
with finite linear combinations of functions of the kind
$g_1(x_{1})\ldots g_r(x_{r})$. By the choice of ${\rm F}_k$ and
the assumptions on $(\gamma_n)$
$$
 {\rm E}^0g_1(\gamma_{k-r+1})\ldots
 g_r(\gamma_{k})=\prod_{i=1}^r{\rm E} g_i(\gamma_{1}).
$$
Hence and from the above assumption on $(\gamma_k)$ condition {\bf
CP2} emerges.
\smallskip

\noindent{\bf Acknowledgement.} The author is grateful to A.
Yurachkivsky for helpful advices.
\smallskip

\begin{center}
{\sc Bibliography}
\end{center}

\begin{list}{\arabic{N}.}{\parsep=-1mm \usecounter{N}}
{\small

\item[1.] Dorogovtsev A. \ Ya., {\it Estimation theory for
parameters of random processes (Russian).} Kyiv University Press.
Kyiv (1982).

\item[2.] Yurachkivsky A. \ P., Ivanenko D. \ O., {\it Matrix
parameter estimation in an autoregression model with
non-stationary noise (Ukranian),} Th. Prob. Math. Stat. 72 (2005),
158--172.

\item[3.] Elsholz, L.\ E., {\it Differential equations and
calculus of variations (Russian),} Nauka, Moscow (1969).

\item[4.]  Chung K. \ L., Williams R. \ J., {\it Introduction in
stochastic integration (Russian),} Mir,  Moscow (1987).

\item[5.]  Yurachkivsky A. \ P., Ivanenko D. \ O., {\it Matrix
parameter estimation in an autoregression model,} Theory of
Stochastic Processes {\bf 12(28)} No 1-2 (2006), 154-161.

\item[6.]  Liptser R. \ Sh.,  Shiryaev A. \ N., {\it Theory of
martingales (Russian),} Nauka,  Moscow (1986).

\item[7.]  Kelley, J., {\it General topology (Russian).}  Nauka,
Moscow (1981).

}
\end{list}

{\small \noindent{\sc Department of Mathematics and Theoretical
Radiophysic,\linebreak Kyiv National Taras Shevchenko University,
Kyiv, Ukraine}

\noindent{\it E-mail:} ida@univ.kiev.ua}

\end{document}